\documentclass[a4paper,11pt]{article}
\usepackage{latexsym}
\usepackage{amssymb}
\usepackage{enumerate}
\usepackage{amsfonts}
\usepackage{amsmath}
\usepackage[dvips]{graphicx}
\usepackage[paper=a4paper,left=30mm,right=20mm,top=25mm,bottom=30mm]{geometry}

\newcommand{\qed}{$\Box$}

\newenvironment{@abssec}[1]{%
    \if@twocolumn

      \section*{#1}%
    \else

      \vspace{.05in}\footnotesize
      \parindent .2in
 {\upshape\bfseries #1. }\ignorespaces
    \fi}

    {\if@twocolumn\else\par\vspace{.1in}\fi}

\newenvironment{keywords}{\begin{@abssec}{\keywordsname}}{\end{@abssec}}

\newenvironment{AMS}{\begin{@abssec}{\AMSname}}{\end{@abssec}}

\newcommand\keywordsname{Key words}
\newcommand\AMSname{AMS subject classifications}
\newcommand\AMname{AMS subject classification}

\def\qed{\vbox{\hrule height0.6pt\hbox{%
  \vrule height1.3ex width0.6pt\hskip0.8ex
  \vrule width0.6pt}\hrule height0.6pt
 }}

  \newtheorem{thm}{Theorem}[section]
  \newtheorem{prop}[thm]{Proposition}
  
  \newtheorem{lem}[thm]{Lemma}

  \newtheorem{rem}[thm]{Remark}

\def\bR{{\mathbf{R}}} 

\def\bN{{\mathbf{N}}}

\title{A simple proof of a strong comparison principle for semicontinuous viscosity solutions of the prescribed mean curvature equation
\thanks{This research was partially supported by the Grants-in-Aid
for Challenging Exploratory Research ($\sharp$ 25610024) and for Scientific Research (B) ($\sharp$ 26287020) of
Japan Society for the Promotion of Science.}}

\author{Masaki {\sc Ohnuma}\thanks{Department of Mathematical Sciences,
Graduate School of Science and Technology, Tokushima University, Tokushima, 770-8502,  Japan. 
({\tt ohnuma@tokushima-u.ac.jp}).}  and Shigeru {\sc Sakaguchi}\thanks{Research Center for Pure and Applied Mathematics,
Graduate School of  Information Sciences, Tohoku
University, Sendai, 980-8579,  Japan.
({\tt sigersak@tohoku.ac.jp}).}}
\begin{document}

\maketitle

\begin{abstract}
A strong comparison principle for semicontinuous viscosity solutions of the prescribed 
mean curvature equation is considered. The difficulties of the problem come from the fact that this nonlinear equation is 
non-uniformly elliptic, does not depend on the value of unknown functions, depends on spatial variables 
and solutions are semicontinuous. 
Our simple proof of the strong comparison principle consists only of three ingredients, the definition 
of viscosity solutions, the inf and sup convolutions of functions, and the theory of classical solutions 
of quasilinear elliptic equations. 
Once we have the strong comparison principle, 
we can prove a weak comparison principle for semicontinuous viscosity solutions of the prescribed 
mean curvature equation in a bounded domain.
 \end{abstract}




\pagestyle{plain}
\thispagestyle{plain}

\section{Introduction}
We consider the prescribed mean curvature equation of the form
\begin{equation}
\label{prescribed mean curvature equation}
\mbox{div}\left(\frac{Du}{\sqrt{1+|Du|^2}}\right)=NH \quad \mbox{ in }\ \Omega, 
\end{equation}
where $\Omega$ is a domain in $\bR^N$ and $N\ge 2$. The function $u:\Omega\to \bR$ is unknown, 
$Du$ denotes the gradient of $u$ in spatial variables $x$ and $H:\Omega\to \bR$ is a locally Lipschitz continuous function in $\Omega$. 
When the solution $u$ is Lipschitz continuous,  equation \eqref{prescribed mean curvature equation} is regarded as uniformly elliptic. 
However, when $u$ is only semicontinuous,  equation \eqref{prescribed mean curvature equation} can be non-uniformly elliptic.

Our goal is to prove a strong comparison principle for semicontinuous viscosity solutions 
of the prescribed mean curvature equation. Here, our strong comparison principle is stated as follows: if 
a lower semicontinuous viscosity supersolution $u$ and an upper semicontinuous viscosity subsolution $v$ satisfy that 
$u\ge v$ in $\Omega$ and $u(x_0)=v(x_0)$ at some point $x_0\in\Omega$,  then $u\equiv v$ in $\Omega$. 

It is well known that for linear elliptic equations the strong comparison principle is equivalent to 
the strong maximum principle since the difference of two solutions is still a solution. Here, the strong 
maximum principle is the following: if a subsolution $u$ satisfies that $u\le m$ with some constant $m$ and 
$u(x_0)=m$ at some point $x_0\in\Omega$, then $u\equiv m$ in $\Omega$. Evidently the strong comparison 
principle implies the strong maximum principle provided that the constant $m$  is a supersolution. The strong 
maximum principle for classical solutions of linear and nonlinear elliptic equations has been well studied 
(cf. \cite{GT}, \cite{PW}). In a book \cite[Theorem 2.1.3 (Tangency Principle),\ p. 16]{PS} we can find 
the strong comparison principle for classical solutions of nonlinear elliptic equations. 

There are some results on the strong maximum principle for weak solutions in the viscosity sense. 
For notations of viscosity solutions we refer to the literature \cite{CIL} and \cite{Ko}. 
The strong maximum principle for semicontinuous viscosity solutions has been proved by \cite{KaKu}, \cite{BD},
\cite{BB}, \cite{KoKo}, and \cite{BGI}.
There are a few papers on the strong comparison principle. Trudinger \cite{T1} proved the strong 
comparison principle for Lipschitz continuous viscosity solutions of uniformly elliptic equations. 
Ishii and Yoshimura \cite{IY} proved the strong comparison principle for semicontinuous 
viscosity solutions of uniformly elliptic equations. At the same time
Giga and the first author \cite{GO} dealt with the strong comparison principle for semicontinuous viscosity 
solutions of nonlinear elliptic equations.  We recently noticed that the argument in \cite[Proof of Theorem 3.1, pp. 177--179]{GO} works for uniformly elliptic equations of the form $F(D^2u)=0$, but  it does not work  for non-uniformly elliptic equations of the form $F(Du, D^2u) = f(x)$ such as \eqref{prescribed mean curvature equation}. 

In the present paper we consider  lower semicontinuous viscosity supersolutions 
and upper semicontinuous viscosity subsolutions of \eqref{prescribed mean curvature equation}. Therefore, we have to deal with non-uniformly elliptic equations.
Our proof is different from usual one. 
After being reduced  to the case where both the supersolution $u$ and the subsolution $v$ are bounded,  by virtue of Jensen, Lions and Souganidis \cite{JLS}, we introduce the inf and sup convolutions of  $u$ and $v$ respectively, where  
 those convolutions are continuous functions and moreover they are monotone with respect to the parameter. 
Then we consider the Dirichlet problems for \eqref{prescribed mean curvature equation} in every sufficiently small ball centered at a point $x_0$, where $u$  
touches $v$.  We choose the continuous boundary data as
the inf and sup convolutions of  $u$  and $v$, respectively. Since $H$ is locally Lipschitz continuous,  by the theory of 
quasilinear elliptic equations (see \cite{GT}),  the gradient estimates of classical solutions are available and these problems have unique classical solutions provided that the ball is sufficiently small.
Here the strong comparison principle is applicable to these two classical solutions and also a weak comparison principle is applicable to compare $u$ and $v$  with these two classical solutions, respectively. Eventually, these comparisons yield that  $u$ and $v$ coincide with each other on the boundary of each small ball centered at a point $x_0$, and hence $u$ and $v$ coincide with each other in a small ball centered at a point $x_0$. Then the conclusion follows from the connectedness of the domain.

The present paper is organized as follows. In section 2 we state our main theorem and prove it. 
In section 3 we give a weak comparison principle as a corollary of our strong comparison principle. 
In Appendix we prove a weak comparison principle for \eqref{prescribed mean curvature equation}
which compares a lower semicontinuous viscosity supersolution with a classical solution, or an upper semicontinuous viscosity subsolution with a classical solution.

\setcounter{equation}{0}
\setcounter{theorem}{0}

\section{Main theorem}
\noindent
Let $\Omega$ be a domain in $\bR^N$, $N\ge 2$ and let $u:\Omega \to \bR$. For functions $u$ we set
$$
M(u):=\mbox{div}\left(\frac{Du}{\sqrt{1+|Du|^2}}\right).
$$
Here $Du$ denotes the gradient of $u$ in spatial variables $x$. 
Let $H:\Omega\to \bR$ be a locally Lipschitz continuous function in $\Omega$. 
Then equation \eqref{prescribed mean curvature equation} is written as
\begin{equation}
\label{PMC}
M(u)=NH \quad \mbox{ in }\ \Omega.
\end{equation}
\par
Our main theorem concerns an extension of the strong comparison theorem to semicontinuous 
viscosity supersolutions and  subsolutions of \eqref{PMC}.
We will use the following notations:
$$
 \begin{array}{ll}
\mbox{USC}(\Omega)&=\{\mbox{upper semicontinuous functions}\ u: \Omega \to \bR \}, \\
\mbox{LSC}(\Omega)&=\{\mbox{lower semicontinuous functions}\ u: \Omega \to \bR \}.
 \end{array}
$$
Also, $\mbox{USC}(\overline{\Omega}),  \mbox{LSC}(\overline{\Omega})$ are defined similarly.

\begin{thm}
\label{main theorem}
Let $u\in \mbox{\rm LSC}(\Omega)$ be a viscosity supersolution of  \eqref{PMC},  that is, 
$$
 M(u) \le NH\quad \mbox{ in } \ \Omega
$$
in the viscosity sense, and let $v\in \mbox{\rm USC}(\Omega)$ be a viscosity subsolution of \eqref{PMC},  that is,
$$
  M(v) \ge NH \quad \mbox{ in } \  \Omega
$$
in the viscosity sense. Assume that $u\ge v$ in $\Omega$ and that $u(x_0)=v(x_0)$ at some point $x_0\in\Omega$. Then $u\equiv v$ in $\Omega$.
\end{thm}
\begin{rem} A continuous viscosity solution $u$ of \eqref{PMC} means that $u \in C(\Omega)$ is both a viscosity supersolution and  subsolution of \eqref{PMC}.
Combining the results of {\rm \cite{B}} and {\rm\cite{T1}} yields the strong comparison principle for 
continuous viscosity solutions of \eqref{PMC}.  Indeed, it is shown in {\rm\cite{B}}  that 
continuous viscosity solutions of \eqref{PMC} are Lipschitz continuous. Then,  equation \eqref{PMC} is regarded as uniformly elliptic, and hence thanks to 
Trudinger's results in  {\rm\cite{T1} } we see that the strong comparison principle for continuous 
viscosity solutions of \eqref{PMC} holds.
\end{rem}
The following weak comparison principle, which is proved in more general form in 
\cite[Theorem 3,\ p. 475]{KaKu},  plays a key role in the present paper. 
Therefore we give a simple proof directly by using the implicit function theorem and the definition of viscosity solutions in the Appendix.
\begin{prop}
\label{prop:a weak comparison by KaKu}
Let $\Omega$ be a bounded domain in $\bR^N$. Let $u\in \mbox{\rm LSC}(\overline{\Omega})$ 
be a viscosity supersolution of \eqref{PMC} and let $v\in C^2(\Omega)\cap C(\overline{\Omega})$ 
be a classical solution of \eqref{PMC}. 
Assume that $u\ge v$ on $\partial\Omega$, then $u\ge v$ in $\Omega$. Similarly, it 
holds for a classical solution $u\in C^2(\Omega)\cap C(\overline{\Omega})$ and 
a viscosity subsolution $v\in \mbox{\rm USC}(\overline{\Omega})$ of \eqref{PMC}.
\end{prop}

Now we are in a position to prove Theorem \ref{main theorem}. 
\bigskip

\noindent
{\it Proof of Theorem \ref{main theorem}.}

\noindent
{\bf 1st step: Reduction to the case where both $u$ 
and $v$ are bounded}: Let $E$ be a bounded domain with $\overline{E}\subset\Omega$ and $x_0\in E$. 
Since $u$ is lower semicontinuous and $v$ is upper semicontinuous, there exists $K>0$ such that
$$
u>-K \  \mbox{ and } \ v<K \quad \mbox{ in }\ \overline{E}.
$$
We use a notation $B(x, r)$ as an open ball in $\bR^N$ of radius $r > 0$ centered at $x\in\bR^N$. 
For simplicity we write in particular $B_r:=B(x_0, r)$ for every $r > 0$.
Choose a positive number $R > 0$ satisfying 
$$
-\frac{1}{R}\le H \le \frac{1}{R} \quad \mbox{ in } \ \overline{E} \quad 
\mbox{ and } \ \overline{B_R}\subset E.
$$
For $x\in\overline{B_R}$ we set
$$
\begin{array}{ll}
 u_R(x)&:=\min\{u(x),\ K+\sqrt{R^2-|x-x_0|^2}\}, \\
 v_R(x)&:=\max\{v(x),\ -K-\sqrt{R^2-|x-x_0|^2}\}.
 \end{array}
$$
Then $u_R\in\mbox{LSC}(\overline{B_R})$ and $v_R\in\mbox{USC}(\overline{B_R})$ are 
bounded in $\overline{B_R}$. Moreover, $u_R$ is a viscosity supersolution of \eqref{PMC} in 
$B_R$ and $v_R$ is a viscosity subsolution of \eqref{PMC} in $B_R$.
\par
Since $v\le u$ in $\Omega$ and $v(x_0)=u(x_0)$, we have $v\le v_R\le u_R\le u$ in 
$\overline{B_R}$ and $v_R(x_0)=u_R(x_0)$. 
By the definition of $u_R$ and $v_R$ we see that if $u_R\equiv v_R$ in $\overline{B_R}$, 
then $u\equiv v$ in $\overline{B_R}$.
Thus we may assume that $u$ and $v$ are bounded in $\overline{B_R}$.
\par
\smallskip

\noindent
{\bf 2nd step: Introducing the inf and sup convolutions of the super and subsolutions}: 
We introduce the inf and sup convolutions of $u$ and $v$, respectively, as in \cite[Proof of Proposition 2,\ p. 977]{JLS}. 
For small $\varepsilon>0$, we set
$$
\begin{array}{ll}
u_{\varepsilon}(x)&:=\inf\limits_{y\in\overline{B_R}}\left\{u(y)+\frac{|x-y|^2}{2\varepsilon}\right\}\quad 
\mbox{ for }\ x\in\overline{B_R}, \\
v^{\varepsilon}(x)&:=\sup\limits_{y\in\overline{B_R}}\left\{v(y)-\frac{|x-y|^2}{2\varepsilon}\right\}\quad 
\mbox{ for }\ x\in\overline{B_R}.
\end{array}
$$
Notice that $u_{\varepsilon},\ v^{\varepsilon}\in C(\overline{B_R})$,  and at each point $x\in\overline{B_R}$ the inf convolution $u_{\varepsilon}(x)$ increases to
$u(x)$ and the sup convolution $v^{\varepsilon}(x)$ decreases to $v(x)$ as $\varepsilon$ decreases to $0$.
\par

\begin{prop}
\label{inequality for approximations}
For each $\varepsilon>0$, $v^{\varepsilon}\ge u_{\varepsilon}$ in $\overline{B_{\frac{R}{2}}}$.
\end{prop}

\noindent
{\it Proof}. By setting $\rho=\frac{R}{2}$, we observe that for every $0 < r \le \rho$
\begin{equation}
\label{sufficient condition for the existence of Classical solutions}
\frac{1}{r}\ge\frac{N}{N-1}|H|\ \mbox{ in } \partial B_r\ \mbox{ and }\ \int_{B_r} |H|^N dx < \omega_N,
\end{equation} 
where $\omega_N$ denotes the Lebesgue measure of the unit ball in $\bR^N$.

Fix $\varepsilon_0>0$ arbitrarily. Let us show that
$$
v^{{\varepsilon}_0}\ge u_{{\varepsilon}_0}\quad \mbox{ in }\ \overline{B_{\rho}}.
$$
For each $\varepsilon\in(0, {\varepsilon}_0]$, we set
$$
{\delta}_{\varepsilon}:=\min_{y\in\overline{B_{\rho}}}(u_{\varepsilon}(y)-v^{\varepsilon}(y)).
$$
Since $u_{\varepsilon}-v^{\varepsilon}$ is continuous in $\overline{B_R}$, ${\delta}_{\varepsilon}$ 
is well defined. By observing that
$$
\min_{y\in\overline{B_{\rho}}}(u_{\varepsilon}(y)-v^{\varepsilon}(y))\le 
u_{\varepsilon}(x_0)-v^{\varepsilon}(x_0)\le u(x_0)-v(x_0)=0,
$$
we know ${\delta}_{\varepsilon}\le 0$. Since
$u_{\varepsilon}(x)$ increases to
$u(x)$ and $v^{\varepsilon}(x)$ decreases to $v(x)$ as $\varepsilon$ decreases to $0$ at each $x\in\overline{B_R}$, 
${\delta}_{\varepsilon}$ is monotone increasing as $\varepsilon$ decreases to $0$.  Let us show a lemma.
\begin{lem}
\label{lem:delta-epsilon tends to zero}
$$
\lim_{\varepsilon\to 0}{\delta}_{\varepsilon}=0.
$$
\end{lem}
\noindent
{\it Proof}.  We may set
$$
\lim\limits_{\varepsilon\to 0}{\delta}_{\varepsilon}=-\lambda
$$
for some number $\lambda\ge 0$. For each $\varepsilon>0$ there exists a point $y_\varepsilon\in \overline{B_{\rho}}$ 
such that ${\delta}_{\varepsilon}=u_{\varepsilon}(y_\varepsilon)-v^{\varepsilon}(y_\varepsilon)$ and moreover there exist 
points $y_{1, \varepsilon}, y_{2, \varepsilon}\in \overline{B_R}$ such that
$$
\begin{array}{ll}
u_{\varepsilon}(y_\varepsilon)&:=u(y_{1, \varepsilon})+\frac{|y_\varepsilon-y_{1, \varepsilon}|^2}
{2\varepsilon}, \\
v^{\varepsilon}(y_\varepsilon)&:=v(y_{2, \varepsilon})-\frac{|y_\varepsilon-y_{2, \varepsilon}|^2}
{2\varepsilon}.
\end{array}
$$
Since $u_{\varepsilon}(y_\varepsilon), v^{\varepsilon}(y_\varepsilon), u(y_{1, \varepsilon})$ and 
$v(y_{2, \varepsilon})$ are bounded, we must have 
\begin{equation}
\label{fundamental convergences}
y_\varepsilon-y_{1, \varepsilon}\to 0\ \mbox{ and } \quad y_\varepsilon-y_{2, \varepsilon}\to 0 \quad \mbox{ as } 
\  \varepsilon\to 0. 
\end{equation}
On the other hand, the Bolzano-Weierstrass theorem yields that there exist a sequence $\{\varepsilon_j\}$ which decreases to $0$ 
as $j\to\infty$ and $x_*\in \overline{B_{\rho}}$ satisfying 
\begin{equation}
\label{convergence of subsequence}
y_{\varepsilon_j}\to x_* \quad \mbox{ as }\ j\to\infty. 
\end{equation}
Then it follows from \eqref{fundamental convergences} and \eqref{convergence of subsequence}  that
\begin{equation}
\label{further convergences}
y_{1, \varepsilon_j},\ y_{2, \varepsilon_j} \to x_* \quad \mbox{ as }\  j\to\infty. 
\end{equation}

Since $\dfrac{|y_{\varepsilon_j}-y_{1, \varepsilon_j}|^2}{2\varepsilon_j}$ and
$\dfrac{|y_{\varepsilon_j}-y_{2, \varepsilon_j}|^2}{2\varepsilon_j}$ are bounded, by taking a subsequence 
if necessary, we may suppose that
$$
\frac{|y_{\varepsilon_j}-y_{1, \varepsilon_j}|^2}{2\varepsilon_j}\to \beta_1(\ge 0) 
\quad \mbox{and} \quad 
\frac{|y_{\varepsilon_j}-y_{2, \varepsilon_j}|^2}{2\varepsilon_j}\to \beta_2(\ge 0) 
\quad \mbox{as} \quad j\to\infty
$$
for some numbers $\beta_1, \beta_2$.
The lower semicontinuity of $u$ and $-v$ at $x_*$ yields that,  for every $\eta>0$,  there exists $\gamma>0$ such that 
if $|x-x_*|<\gamma$ then
$$
u(x)>u(x_*)-\eta \quad \mbox{and} \quad -v(x)>-v(x_*)-\eta.
$$
By \eqref{further convergences} there exists $n_0\in \bN$ such that if $j\ge n_0$ then
$$
u(y_{1, \varepsilon_j})>u(x_*)-\eta \ \mbox{ and } \  -v(y_{2, \varepsilon_j})>-v(x_*)-\eta.
$$
 Hence, for $j\ge n_0$
$$
\begin{array}{ll}
\delta_{\varepsilon_j}
&=u_{\varepsilon_j}(y_{\varepsilon_j})-v^{\varepsilon_j}(y_{\varepsilon_j})
=u(y_{1, \varepsilon_j})+\frac{|y_{\varepsilon_j}-y_{1, \varepsilon_j}|^2}{2\varepsilon_j}
-v(y_{2, \varepsilon_j})+\frac{|y_{\varepsilon_j}-y_{2, \varepsilon_j}|^2}{2\varepsilon_j} \\
&>u(x_*)-v(x_*)-2\eta
+\frac{|y_{\varepsilon_j}-y_{1, \varepsilon_j}|^2}{2\varepsilon_j}
+\frac{|y_{\varepsilon_j}-y_{2, \varepsilon_j}|^2}{2\varepsilon_j}.
\end{array}
$$
Letting $j\to\infty$ yields that for every $\eta>0$
$$
0\ge-\lambda\ge u(x_*)-v(x_*)-2\eta+\beta_1+\beta_2.
$$

Since $u(x_*)\ge v(x_*),\ \beta_1\ge 0$ and $\beta_2\ge 0$, we see that $0\le\lambda\le 2\eta$ and 
$0\le\beta_1+\beta_2\le 2\eta$ for every $\eta>0$. Thus we conclude that  $\lambda=\beta_1=\beta_2=0$ and the proof of Lemma \ref{lem:delta-epsilon tends to zero} is finished. 
\qed

\vskip 2ex
\par
We return to the proof of Proposition \ref{inequality for approximations}. 
In order to prove that 
$$
v^{{\varepsilon}_0}\ge u_{{\varepsilon}_0}\  \mbox{ in }\  \overline{B_{\rho}},
$$
 let us show that  for each $r$ with $0<r\le \rho$, 
$$
v^{\varepsilon_0}(x)\ge u_{\varepsilon_0}(x) \  \mbox{ for all } \  x\in\partial B_r.
$$

Suppose that there exist  $r$ with $0<r\le \rho$ and a point $x_1\in \partial B_r$ so that 
$u_{{\varepsilon}_0}(x_1)-v^{{\varepsilon}_0}(x_1)>0$. Since $u_{{\varepsilon}_0}(x)-v^{{\varepsilon}_0}(x)$ 
is continuous, there exist $r_1$ with $0<r_1<\frac{r}{4}$ and $\beta > 0$ such that 
\begin{equation}
\label{strict positivity beta}
u_{{\varepsilon}_0}(x)-v^{{\varepsilon}_0}(x)\ge \beta\ \mbox{  in } \overline{B(x_1, r_1)}.
\end{equation}
We divide $\partial B_r$ into two  pieces:
$$
\Gamma_+:=\partial B_r\cap \overline{B(x_1, r_1)}\ \mbox{ and }\  \Gamma_0:=\partial B_r\setminus\Gamma_+.
$$
Clearly $\partial B_r$ is a disjoint union of $\Gamma_+$ and $\Gamma_0$.
\par
Note that  \eqref{strict positivity beta} gives in particular
\begin{equation}
\label{strict positivity beta 2}
u_{{\varepsilon}_0}(x)\ge v^{{\varepsilon}_0}(x)+\beta \quad \mbox{ on } \  \Gamma_+ .
\end{equation} 
Then for every $0<\varepsilon\le\varepsilon_0$, by the monotonicity of $u_{\varepsilon}$ and $v^{\varepsilon}$
$$
u_{\varepsilon}(x)\ge u_{{\varepsilon}_0}(x)\ge v^{{\varepsilon}_0}(x)+\beta\ge v^{\varepsilon}(x)+\beta
\quad \mbox{ on } \ \Gamma_+,
$$
and by the definition of $\delta_{\varepsilon}$ we see that
\begin{equation}
\label{key approximation inequality}
u_{\varepsilon}(x)\ge v^{\varepsilon}(x)+\delta_{\varepsilon} \quad \mbox{ on } \ \partial B_r. 
\end{equation}

Since $H$ is Lipschitz continuous in $\overline{B_r}$,  we see that the interior estimates of \cite[Corollary 16.7,  p. 407]{GT}  are available. Therefore it follows from \eqref{sufficient condition for the existence of Classical solutions} and the theory of 
the prescribed mean curvature equation (\cite[Theorem 16.10, p. 408 ]{GT}) that 
there exist $\hat u_{\varepsilon},\ \hat v^{\varepsilon}\in C^{2}(B_r)\cap C(\overline{B_r})$ 
satisfying
$$
\begin{array}{ll}
M(\hat u_{\varepsilon})=M(\hat v^{\varepsilon})=NH \quad &\mbox{ in } \ B_r, \\
\hat u_{\varepsilon}=u_{\varepsilon} \quad \mbox{and} \quad \hat v^{\varepsilon}=v^{\varepsilon} 
\quad &\mbox{ on } \ \partial B_r.
\end{array}
$$
Notice that $M(\hat v^{\varepsilon}+\delta_{\varepsilon})=NH$ in $B_r$. By Proposition \ref{prop:a weak comparison by KaKu}
and \eqref{key approximation inequality} we observe that
\begin{equation}
\label{inequalities including delta-epsilon}
v+\delta_{\varepsilon}\le \hat v^{\varepsilon}+\delta_{\varepsilon}\le \hat u_{\varepsilon}\le u 
\quad \mbox{ in } \ \overline{B_r}. 
\end{equation}
Also, $\hat u_{\varepsilon}$ increases and $\hat v^{\varepsilon}$ decreases  
as $\varepsilon$ decreases to $0$ by the monotonicity of $u_{\varepsilon}$ and $v^{\varepsilon}$.

The boundedness of $\{\hat u_{\varepsilon}\}$ and $\{\hat v^{\varepsilon}\}$ together with the interior estimates of \cite[Corollary 16.7, p. 407]{GT} 
yields that there exist 
$\hat u_0, \hat v^0\in C^{2}(B_r)$ such that
$$
\begin{array}{ll}
&\hat u_{\varepsilon}\to \hat u_0 \  \mbox{ and } \  \hat v^{\varepsilon}\to \hat v^0 \  \mbox{ as } \  \varepsilon \to 0 \mbox{ uniformly on compact sets in } \ B_r, 
\\
&M(\hat u_0)=M(\hat v^0)=NH 
\ \mbox{ in } \ B_r.
\end{array}
$$

Therefore, since $\lim\limits_{\varepsilon\to 0}\delta_{\varepsilon}=0$,  we observe from \eqref{inequalities including delta-epsilon} that
$$
v\le \hat v^0\le \hat u_0\le u 
\quad \mbox{ in } \ B_r.
$$
Since $v(x_0)=u(x_0)$, we have that $\hat v^0(x_0)=\hat u_0(x_0)$. By the strong comparison principle for classical solutions 
we see that
\begin{equation}
\label{the equal limits}
\hat v^0(x)\equiv \hat u_0(x)  \quad \mbox{in} \  B_r. 
\end{equation}

However,  by \eqref{strict positivity beta 2}
$$
\hat u_{{\varepsilon}_0}(x)\ge\hat v^{{\varepsilon}_0}(x)+\beta \quad \mbox{ on } \ \Gamma_+.
$$
Hence it follows from  the continuity of $\hat u_{{\varepsilon}_0}$ and $\hat v^{{\varepsilon}_0}$ that there exist $\beta_3$ with 
$0<\beta_3\le\beta$  and $r_2$ with $0 < r_2 \le r_1$ satisfying
$$
\hat u_{{\varepsilon}_0}(x)\ge\hat v^{{\varepsilon}_0}(x)+\beta_3 
\quad \mbox{ in } \  B_r\cap \overline{B(x_1,r_2)}.
$$
By the monotonicity of $u_{\varepsilon}$ and $v^{\varepsilon}$ for $0<\varepsilon\le\varepsilon_0$
we observe that
$$
\hat u_{\varepsilon}(x)\ge\hat u_{{\varepsilon}_0}(x)\ge\hat v^{{\varepsilon}_0}(x)+\beta_3
\ge\hat v^{\varepsilon}(x)+\beta_3
\quad \mbox{ in } \ B_r\cap \overline{B(x_1,r_2)}.
$$
Therefore letting $\varepsilon\to 0$ yields that
$$
\hat u_0(x)\ge\hat v^0(x)+\beta_3 
\quad \mbox{ in } \ B_r\cap \overline{B(x_1,r_2)}.
$$
This contradicts \eqref{the equal limits}
.
\par
Eventually, we conclude that  $v^{\varepsilon_0}\ge u_{\varepsilon_0}$ on $\partial B_r$. 
Since this holds for every $0<r\le \rho$, we have
$$
v^{\varepsilon_0}\ge u_{\varepsilon_0} 
\quad \mbox{ in } \ \overline{B_{\rho}}. \quad \mbox{\qed}
$$

\vskip 2ex
\noindent
{\bf 3rd step: Completion of the proof of Theorem \ref{main theorem}}: 
By Proposition \ref{inequality for approximations},  letting $\varepsilon\to 0$ yields that
$$
v\ge u 
\quad \mbox{ in } \  \overline{B_{\rho}},
$$
which shows that $u\equiv v$ in $\overline{B_{\rho}}$. Since $\Omega$ is connected, 
we conclude that $u\equiv v$ in $\Omega$. \qed

\setcounter{equation}{0}
\setcounter{theorem}{0}

\section{A weak comparison principle for semicontinuous viscosity solutions of the prescribed 
mean curvature equation}
\noindent
 The strong comparison principle proved in section 2 yields a weak comparison principle for 
semicontinuous viscosity solutions of the prescribed mean curvature equation in a bounded domain.
\begin{thm}
Let $\Omega$ be a bounded domain in $\bR^n$. 
Let $u\in \mbox{\rm LSC}(\overline{\Omega})$ and $v\in \mbox{\rm USC}(\overline{\Omega})$ be viscosity 
super and subsolutions of \eqref{PMC}, respectively. 
Assume that $u\ge v$ on $\partial\Omega$. Then $u\ge v$ in $\Omega$, 
and hence either $u\equiv v$ in $\Omega$ or $u>v$ in $\Omega$.
\end{thm}

\noindent
{\it Proof.}\quad 
Suppose that there exists a point $x_1\in\overline{\Omega}$ satisfying
$$
\theta:=\min_{x\in\overline{\Omega}}(u-v)(x)=(u-v)(x_1)<0.
$$
Hence $x_1\in\Omega$, since $u\ge v$ on $\partial\Omega$. Then we observe that
$$
u\ge v+\theta \quad \mbox{ in }\  \Omega \quad \mbox{and}\quad 
u(x_1)=v(x_1)+\theta.
$$
Note that $v+\theta$ is also a viscosity subsolution of \eqref{PMC}. By Theorem  \ref{main theorem} we have that $u\equiv v+\theta$ in $\Omega$, which contradicts
the assumption that $u\ge v$ on $\partial\Omega$. Therefore we see that $u\ge v$ in $\Omega$.
\par
Moreover, if there exists a point $x_2\in\Omega$ so that $u(x_2)=v(x_2)$, then $u\equiv v$ in $\Omega$ 
by Theorem \ref{main theorem}, which concludes that either $u\equiv v$ in $\Omega$ or 
$u>v$ in $\Omega$. \qed




\renewcommand{\thethm}{\Alph{thm}} 

\section*{Appendix}
\noindent
Although  Proposition \ref{prop:a weak comparison by KaKu} was already proved in \cite[Theorem 3,\ p. 457]{KaKu}, 
for convenience we will give a simple proof directly by using the implicit function theorem and the definition of 
viscosity solutions.

\noindent
{\it Proof of Proposition \ref{prop:a weak comparison by KaKu}.}\quad By the argument in \cite[Theorem A.1,\ p. 253]{MS}, which applies  Sard's theorem to a smooth function being comparable to the distance function to the closed set $\mathbb R^N \setminus \Omega$ due to Calder\'on and Zygmund \cite[Lemma 3.6.1,
p. 136]{Z} (see also \cite[Lemma 3.2, p. 185]{CZ}),  we observe that 
for each small 
$\varepsilon>0$ there exists a smooth open set $\Omega_{\varepsilon}\subset\subset\Omega$ with 
$\Omega_{\varepsilon'}\subset\Omega_{\varepsilon}$ if $\varepsilon<\varepsilon'$ and 
$\Omega_{\varepsilon}\to\Omega$ as $\varepsilon\to 0$. Since $\Omega$ is bounded, we notice that $\Omega_{\varepsilon}$ is 
a union of a finite number of smooth domains. 

Let $u\in \mbox{\rm LSC}(\overline{\Omega})$ be a viscosity supersolution of \eqref{PMC} and let $v\in C^2(\Omega)\cap C(\overline{\Omega})$ be a classical solution of \eqref{PMC}. 
Assume that $u\ge v$ on $\partial\Omega$.  Since $H$ is locally Lipschitz continuous,  we see that the interior estimates of \cite[Corollary 16.7,  p. 407]{GT}  are available. Therefore,  with the aid of the Schauder interior estimates for elliptic equations,  there exists a number $\alpha$ with $0<\alpha<1$ depending on $\varepsilon$ such that $v\in C^{2,\alpha}(\overline{\Omega_{\varepsilon}})$.  We can write 
\begin{equation}
\label{approximation of the domain}
\Omega_{\varepsilon}=\bigcup_{j=1}^{n(\varepsilon)}D_{\varepsilon, j}.
\end{equation}
Consider an arbitrary $D_{\varepsilon, j}$. For simplicity we will 
write $D$ instead of $D_{\varepsilon, j}$. 
Note that $\partial D$ is close to $\partial \Omega$ if $\varepsilon>0$ is sufficiently small. 
Since $u\in \mbox{\rm LSC}(\overline{\Omega})$, $v\in C(\overline{\Omega})$ and $u\ge v$ on $\partial\Omega$, 
there exists $\tau(\varepsilon)>0$ satisfying $\lim\limits_{\varepsilon\to 0}\tau(\varepsilon)=0$ 
and $u>v-\tau(\varepsilon)$ on $\partial D$.

Set $w_{\varepsilon}:=v-\tau(\varepsilon)$ and we have 
$$
M(w_{\varepsilon})=NH 
\quad \mbox{ in } \ D.
$$
We set
$$
\begin{array}{ll}
X&:=\{f\in C^{2,\alpha}(\overline{D})\ |\ f=0 \ \mbox{on} \ \partial D\}, \\
F&:X\times \bR \ni (f,s)\mapsto M(w_{\varepsilon}+f)-N(H+s)
\in C^{\alpha}(\overline{D}).
\end{array}
$$
We use the implicit function theorem for $X$ and $F$( see \cite[Theorem 2.3,\ p. 38]{AP}, 
\cite[Theorem 15.1,\ p. 148]{D} for instance ). 
For each $0<\delta<<1$ there exists 
$\tilde w_{\varepsilon,\delta}\in C^{2,\alpha}(\overline{D})$ satisfying
$$
\begin{array}{ll}
M(\tilde w_{\varepsilon,\delta})&=N(H+\delta)
\quad \mbox{ in } \  D, \\
\tilde w_{\varepsilon,\delta}&=w_{\varepsilon}
\quad \mbox{ on } \ \partial D.
\end{array}
$$
Note that $\tilde w_{\varepsilon,0}=w_{\varepsilon}$. Then we have
\begin{equation}
\label{conclusion by viscosity supersolution}
u\ge\tilde w_{\varepsilon,\delta}
\quad \mbox{ in } \ \overline{D}. 
\end{equation}
Indeed, suppose that there exists a point $z\in\overline{D}$ satisfying
\begin{equation}
\label{the test function}
\min_{x\in\overline{D}}(u-\tilde w_{\varepsilon,\delta})(x)=(u-\tilde w_{\varepsilon,\delta})(z) <0.
\end{equation}
Hence $z\in D$, since $u-\tilde w_{\varepsilon,\delta}=u-w_{\varepsilon}>0$ on $\partial D$. 
Moreover, since $u$ is a viscosity supersolution of \eqref{PMC}, we have from \eqref{the test function} that 
$M(\tilde w_{\varepsilon,\delta})(z)\le NH(z)$.  This contradicts
the fact that $M(\tilde w_{\varepsilon,\delta})(z)= N(H(z)+\delta)$ with $\delta > 0$.

Letting $\delta\to 0$ in \eqref{conclusion by viscosity supersolution} yields that
$$
u\ge w_{\varepsilon} \quad \mbox{ in } \ \overline{D}.
$$
Hence it follows from \eqref{approximation of the domain} that 
$$
u\ge v-\tau(\varepsilon) \quad \mbox{ in } \  \Omega_{\varepsilon}.
$$
Thus, letting $\varepsilon\to 0$ yields that $u\ge v$ in $\Omega$, which completes the proof of Proposition \ref{prop:a weak comparison by KaKu}.

\vskip 4ex

\section*{Acknowledgements}
The  authors are grateful  to Professor Hiroyoshi Mitake for informing them of a gap in the proof in \cite{GO} when one deals with the prescribed mean curvature equation.

\medskip

\end{document}